%% file: weil1.tex
\documentclass[11pt,twoside]{article}

\usepackage{latexsym}
\usepackage{amsmath}
\usepackage{amsthm}
\usepackage{amssymb}
\usepackage{vmargin}
\usepackage{amscd}
\usepackage{stmaryrd}
\usepackage{euscript}
\usepackage{mathrsfs}
\usepackage{amscd}
\usepackage[all]{xy}
\usepackage{xr}
\DeclareMathAlphabet{\mathpzc}{OT1}{pzc}{m}{it}

\setmargins{32mm}{20mm}{14.6cm}{22cm}{1cm}{1cm}{1cm}{1cm}

\include{newcommands}

\include{newthms}

\begin{document}
\title{Weight decompositions on \'etale fundamental groups}
\author{J.P.Pridham\thanks{The author is supported by Trinity College, Cambridge.}}
\maketitle

\begin{abstract}
Let $X_0$ be a smooth or proper variety defined over a finite field $k$, and let $X$ be the base extension of $X_0$ to the algebraic closure $\bar{k}$. The geometric \'etale fundamental group $\pi_1(X,\bar{x})$ of $X$ is a normal subgroup of the Weil group, so conjugation gives it a Weil action. For $l$ not dividing the characteristic of $k$, we consider the pro-$\Q_l$-algebraic completion of $\pi_1(X,\bar{x})$ as a non-abelian Weil representation. Lafforgue's Theorem and Deligne's Weil II theorems imply that this affine group scheme is mixed, in the sense that its structure sheaf is a mixed Weil representation. When $X$ is smooth, weight restrictions apply, affecting the possibilities for the  structure of this group. This gives new examples of groups which cannot arise as \'etale fundamental groups of smooth varieties. 
\end{abstract}

\tableofcontents

\section*{Introduction}
\addcontentsline{toc}{section}{Introduction}

Let $X_0$  be a connected variety defined over a finite field $k=\bF_q$,  equipped with a point  $x \in X_0(k)$, and let $l$ be a prime not dividing $q$. If we set $X:=X_0\ten_k\bar{k}$, the embedding $i: \pi_1(X,\bar{x}) \into W(X_0,x)$ of the fundamental group into the Weil group gives a conjugation action of the Weil group on the fundamental group. The Weil conjecture for $\H^1(X,\Q_l)$ can be re-expressed by considering $\H^1(X,\Q_l)^{\vee}$ as the $l$-adic Weil representation $V$, universal among continuous  Weil-equivariant group homomorphisms
$$
\pi_1(X,x) \to V,
$$
and then stating that this Weil representation is mixed.

We consider a non-abelian version of this, by defining the pro-$\Q_l$-algebraic group ${}^W\!\algpia$ to be the 
 universal object classifying continuous $W(X_0,x)$-equivariant homomorphisms
$$
\pi_1(X,\bar{x}) \to G(\Q_l),
$$  
where $G$ ranges over all algebraic groups $G$ over $\Q_l$ equipped with continuous $W(X_0,x)$-actions. Representations of $\w\algpia$ are precisely $\pi_1(X,\bar{x})$-subrepresentations of $W(X_0,x)$-representations.

We say that an algebraic Weil action on a pro-algebraic group $G$ is mixed if the structure sheaf  $O(G)$ is a sum of mixed Weil representations. 

The Levi decomposition for pro-algebraic groups allows us to write
$$
{}^W\!\algpia \cong \wmypia\rtimes \w\redpia,
$$
where  $\wmypia$ is the pro-unipotent radical of ${}^W\!\algpia$ and ${}^W\!\redpia$ is the pro-reductive completion of $\w\algpia$. This decomposition is  unique up to conjugation by $\wmypia$. 

In Section \ref{del}, we use Lafforgue's Theorem to show that for any variety,  the Weil action on $\w\redpia$ is pure of weight zero. Deligne's Weil II theorems then show that, if $X$ is smooth or proper,  the Weil action on ${}^W\!\algpia$  is mixed. This can be thought of as a direct analogue of the non-abelian Hodge theorems of \cite{Simpson}. One consequence is that for any morphism $f:X \to Y$ of varieties over $\bF_q$,  with $X$ smooth and $\bV$ any semisimple constructible $\Q_l$-local system underlying a Weil sheaf on $Y$, the pullback $f^{-1}\bV$ is semisimple.

The rest of the paper is  dedicated to studying the Weil action on $\wmypia$ when $X$ is smooth or proper, and thus establishing restrictions on the structure of the fundamental group.
In order to study the pro-unipotent extension $\w\algpia \to \w\redpia$, we use deformation-theoretic machinery. The  group $\wmypia$ is the universal deformation
$$
\rho: \pi_1(X,\bar{x}) \to U \rtimes {}^W\!\redpia
$$
of the canonical representation
$$
\rho_0:\pi_1(X,\bar{x}) \to {}^W\!\redpia,
$$
for $U$ pro-unipotent. In \cite{higgs},  a theory of deformations over nilpotent Lie algebras with $G$-actions was developed,  and this enables us to analyse our scenario. 

In  Section \ref{final}, we use Deligne's Weil II theorems  to study $\wmypia$.
If $X$ is smooth and proper, then the weight decomposition on $\wmypia$ splits the lower central series filtration, and it  is quadratically presented, in the sense that its Lie algebra can be defined by equations of bracket length two. If $X$ is merely smooth, then $\wmypia$ is defined by equations of bracket length at most four. Since rigid representations of the fundamental group extend to Weil representations, these properties are used to give new examples of groups which cannot occur as fundamental groups of smooth varieties in finite characteristic. 

This   generalises the results of \cite{paper1} on deforming reductive representations of the fundamental group. In this paper, we are taking a reductive representation
$$
\rho_0:\pi_1(X,\bar{x}) \to G,
$$
and considering deformations
$$
\rho: \pi_1(X,\bar{x}) \to U \rtimes G
$$
of $\rho_0$, for $U$ unipotent. Effectively, \cite{paper1}  considers only $U=\exp(\mathrm{Lie}(G) \ten \m_A)$, for  $\m_A$ a maximal ideal of an Artinian local $\Q_l$-algebra. Since taking $U=\mypia$ pro-represents this functor when $G=\redpia$, all examples can be understood in terms of the structure of $\mypia$.

The structure result in the smooth and proper case is much the same as those established in \cite{malcev} and \cite{higgs} for fundamental groups of compact K\"ahler manifolds. Likewise, \cite{paper1}  was the  analogue in finite characteristic of Goldman and Millson's results on K\"ahler representations (\cite{gm}).

\section{The pro-algebraic fundamental group as a Weil representation}\label{del}
\subsection{Algebraic actions}
All pro-algebraic groups in this paper will be defined over fields of characteristic zero (usually $\Q_l$). All representations of pro-algebraic groups will be finite-dimensional.

\begin{definition}
Given a pro-algebraic group $G$, let $O(G)$ denote global sections of the structure sheaf of $G$. This is a sum of  $G\by G$-representations, the actions corresponding to right and left translation. Let $E(G)$ be the dual of $O(G)$ --- this is a pro-$G\by G$-representation. In fact, since any coalgebra is the sum of its finite-dimensional subcoalgebras, $E(G)$ is an inverse limit of finite-dimensional (non-commutative) algebras. 

$E(G)$-modules then correspond to pro-$G$-representations, and for a morphism $G \to H$ and a pro-$G$-representation $V$, we define 
$$
\Ind_{G}^{H}V:=V\hat{\ten}_{E(G)}E(H). 
$$
\end{definition}

\begin{definition}
Given a discrete group $\Gamma$ acting on a pro-algebraic group $G$, we define ${}^{\Gamma}\!G$ to be the maximal quotient of $G$ on which $\Gamma$ acts algebraically. This is the inverse limit $\Lim_{\alpha} G_{\alpha}$ over those surjective maps
$$
G \to G_{\alpha},
$$
with $G_{\alpha}$ algebraic, for which the $\Gamma$-action descends to $G_{\alpha}$. 
\end{definition}

\begin{lemma}\label{repchar}
The representations of ${}^{\Gamma}\!G$ are precisely those $G$-representations which arise as $G$-subrepresentations of (finite-dimensional) $G\ltimes \Gamma$-representations. 
\end{lemma}
\begin{proof}
Given  ${}^{\Gamma}\!G \xra{\theta} \GL(V)$, there must exist an algebraic quotient group $G_{\alpha}$ of $G$ to which $\Gamma$ descends, with $\theta$ factoring as ${}^{\Gamma}\!G \to G_{\alpha} \to \GL(V)$. Now, since $G_{\alpha}$ is an algebraic group, $\Aut(G_{\alpha})$ is also, and there is a homomorphism  $G\ltimes \Gamma \to H_{\alpha} :=G_{\alpha} \ltimes \Aut(G_{\alpha})$. Since $G_{\alpha}\into H_{\alpha}$, the $G_{\alpha}$-representation $V$ is a subrepresentation of the pro-$H_{\alpha}$-representation $\Ind_{G_{\alpha}}^{H_{\alpha}}V $, so for some quotient representation $\Ind_{G_{\alpha}}^{H_{\alpha}}V \to W$, the composition $V \to W$ must be injective. Thus $V$ is a subrepresentation of the  $G\ltimes \Gamma$-representation $W$.
 
Conversely, Let $V\le W$ be $G$-representations, with  $W$ a $G\ltimes \Gamma$-representation. If we let $G_{\alpha}$ be the image of $G \to \GL(W)$, then the adjoint action of $\Gamma$ on $\GL(W)$ restricts to an action on $G_{\alpha}$. Since the action of $G$ on $W$ preserves $V$, there is an algebraic map $G_{\alpha} \to \GL(V)$, as required.  
\end{proof}

\begin{definition}
Given a pro-algebraic group $G$, we will denote its reductive quotient by $G^{\red}$; this is the universal object among quotients $G \to H$, with $H$  reductive algebraic. Representations of $G^{\red}$ correspond to semisimple representations of $G$. We write $\Ru(G)$ for the kernel of $G \to G^{\red}$ --- this is called the pro-unipotent radical of $G$.
\end{definition}

\begin{lemma}
${}^{\Gamma}(G^{\red})=({}^{\Gamma}G)^{\red}$. We will hence denote this group by ${}^{\Gamma}G^{\red}$.
\end{lemma}
\begin{proof}
Note that in both cases, representations correspond to those semisimple $G$-representations which arise as $G$-subrepresentations of (finite-dimensional) $G\ltimes \Gamma$-representations.
\end{proof}

The Levi decomposition, proved in \cite{Levi}, states that for every pro-algebraic group $G$, the surjection $G \to G^{\red}$ has a section, unique up to conjugation by $\Ru(G)$, inducing an isomorphism $G\cong \Ru(G)\rtimes G^{\red}$. 

\begin{lemma}\label{innerworks}
Given a pro-algebraic group $G$, an automorphism $F$ of $G$, and an element $g \in G$, the action of $F$ on $G$ is algebraic if and only if the action of $\ad_g \circ F$ is algebraic.
\end{lemma}
\begin{proof}
First note that we have an isomorphism from $G\rtimes\langle \ad_g \circ F \rangle$ to $G\rtimes\langle F \rangle$ fixing $G$, given by sending $\ad_g \circ F$ to $g\cdot F$. Hence, by Lemma \ref{repchar},  ${}^{F}\!G={}^{\ad_g \circ F}\!G$.
\end{proof}

\begin{corollary}\label{algboth}
The action of $F$ on $G$ is algebraic if and only if the corresponding actions on $G^{\red}$ and $\Ru(G)$ are.
\end{corollary}
\begin{proof}
Without loss of generality, by the previous lemma, we may assume that $F$ must preserve the Levi decomposition (following conjugation by a suitable element of $\Ru(G)$). Write $F=F^{\red}F^u$, for $F^u:\Ru(G) \to \Ru(G)$, and $F^{\red}:G^{\red} \to G^{\red}$. By Lemma \ref{repchar} and Tannakian duality, ${}^{F}\!G$ is the image of $G \to (G\rtimes \langle F \rangle)^{\alg}$, the latter group being the pro-algebraic completion of $G\rtimes \langle F \rangle$.

Then note that we have an embedding
$$
(G\rtimes\langle F \rangle)^{\alg} \into (\Ru(G)\rtimes\langle F^u \rangle)^{\alg}  \rtimes (G^{\red}\rtimes\langle F^{\red} \rangle)^{\alg}, 
$$
so the map from $G$ to the group on the left is an embedding if and only if the maps from $G^{\red},\Ru(G)$ to the groups on the right are embeddings.
\end{proof}

\begin{lemma}\label{algchar} Let $F$ act on $G\ltimes U$, for $G$ reductive and $U$ pro-unipotent, with $F$ preserving and acting algebraically on $G$. If we also assume that $\Hom_G(V,U/[U,U])$ is finite-dimensional for all $G$-representations $V$,   then $F$ acts algebraically on $G\ltimes U$. 
\end{lemma}
\begin{proof}

By the previous lemma, it suffices to show that $F$ acts algebraically on $U$. Let $S$ be the set of isomorphism classes of irreducible representations of $G$. Since $F$ acts algebraically on $G$, the $F$-orbits in $S$ are all finite. Let $\fu:=\Lie(U)$, and take the canonical decomposition $\fu=\prod_{s\in S}\fu_s$ of $\fu$ as a $G$-representation. Let $T=S/F$ be the set of $F$-orbits in $S$, giving a weaker decomposition $\fu=\prod_{t\in T}\fu_t$, where $\fu_t =\prod_{s \in T}\fu_s  $. 

$F$ is then an automorphism of $\fu$ respecting this decomposition; let $H$ be the group of all such automorphisms. We then have an embedding
$$
U\rtimes\langle F \rangle \into U \rtimes H,
$$
so it suffices to show that the group $H$ is pro-algebraic, since this embedding must then factor through  $(U\rtimes\langle F \rangle)^{\alg}$.

Choose a $G$-equivariant section to the map $\fu \to \fu/[\fu,\fu]$, and let its image be $V$. The group $H$ is a closed subspace of the space of all linear maps $\Hom_T(V,\fu)$ preserving the $T$-decomposition. The hypothesis implies that $V_s$ is finite-dimensional for all $s \in S$, so  $V_t$ must be finite-dimensional for all $t \in T$, the $F$-orbits being finite. Thus $H$ is an affine group scheme, i.e. a pro-algebraic group, as required.
\end{proof}

\begin{lemma}\label{dual} If $G$ is a pro-algebraic group, and we regard $O(G)$ as a sum of $G$-representations via the left action, then for any $G$-representation $V$, $V^{\vee}\cong \Hom_G(V, O(G))$, with the $G$-action on $V^{\vee}$ coming from the right action on $O(G)$.
\end{lemma}
\begin{proof}
This follows immediately from  \cite{tannaka} II Proposition 2.2, which states that $G$-representations correspond to $O(G)$-comodules. Under this correspondence, $\alpha \in V^{\vee}$ is associated to the morphism which sends $v \in V$ to the function $g \mapsto \alpha(g\cdot v)$. 
\end{proof}

\begin{lemma}\label{fdual}
If an endomorphism $F$ acts on a pro-algebraic group $G$ and compatibly on a $G$-representation $V$ (i.e. $F(g\cdot v)=(Fg)\cdot (Fv)$), then the  dual action of $F$ on $V^{\vee}$ corresponds to the action on $\Hom_G(V, O(G))$ which sends $\theta$ to the composition
$$
V \xra{F} V \xra{\theta} O(G) \xra{F^*} O(G)
$$
 \end{lemma}

\subsection{Weil actions}\label{weilact}

Let $k=\bF_q$, take a  connected variety $X_0/k$, and let $X=X_0\ten_k\bar{k}$. Fix a closed point $x$ of $X$, and denote the associated geometric point $x\ten_{k(x)}\bar{k} \to X$ by $\bar{x}$. Without loss of generality (increasing $q$ if necessary), we assume that $k(x) \subset \bF_q$. Let $l$ be a prime not dividing $q$, and consider the pro-$\Q_l$-algebraic completion $\algpia$ of the  \'etale fundamental group $\pi_1(X,\bar{x})$ of $X$. This is the universal object classifying continuous homomorphisms
$$
\pi_1(X,\bar{x}) \to G(\Q_l),
$$  
where $G$ ranges over all algebraic groups $G$ over $\Q_l$.

Recall that the Frobenius element gives a canonical generator of $\pi_1(\Spec k)\cong \hat{\Z}$, and that the Weil group $W(X_0,x)$ is defined by
$$
W(X_0,x) = \pi_1(X_0,\bar{x})\by_{\hat{\Z}}\Z,
$$
which has $\pi_1(X,\bar{x})$ as a normal subgroup.
Observe that the conjugation action of $W(X_0,x)$ on $\pi_1(X,\bar{x})$ then extends by universality to an  action of $W(X_0,x)$ on $\algpia$. Let $F_x \in W(X_0,x)$ be the Frobenius element associated to $x$.

\begin{lemma}\label{wf}
If $W:=W(X_0,x)$ and $F:=F_x$, then $\w\algpia={}^{F}\!\algpia$, with representations of this group being those continuous $\pi_1(X,\bar{x})$-representations which arise as $\pi_1(X,\bar{x})$-subrepresentations of Weil representations.
\end{lemma}
\begin{proof}
By Lemma \ref{repchar}, representations of $\w\algpia$ are continuous $\pi_1(X,\bar{x})$-subrepresentations of $\pi_1(X,\bar{x})\ltimes W(X_0,x)$-representations. These are precisely $\pi_1(X,\bar{x})$-subrepresentations of $W(X_0,x)$-representations. Since  $W(X_0,x)=\pi_1(X,\bar{x})\ltimes \langle F_x\rangle$, these are the same as representations of ${}^{F}\!\algpia$. By Tannakian duality (\cite{tannaka}), this determines the quotient groups $\w\algpia,{}^{F}\!\algpia$ of $\algpia$, which must then be equal.
\end{proof}

\begin{lemma}
$\w\algpia$ is the image of the homomorphism $\algpia\xra{i} \mathpzc{W}(X_0,x)$, where $\mathpzc{W}(X_0,x)$ is the pro-algebraic completion of the Weil group $W(X_0,x)$. 
\end{lemma}
\begin{proof}
Representations of $\im(i)$ are those $\algpia$ representations $V$ for which $V \to \Ind_{\algpia}^{\mathpzc{W}(X_0,x)}$ is injective. By Lemmas \ref{repchar} and \ref{wf}, these are the same as $\w\algpia$-representations.
\end{proof}

\begin{definition}
Given a pro-$\Q_l$-algebraic group $G$, equipped with an algebraic action of the Weil group $W(X_0,x)$, we will say that this Weil action  on $G$ is  mixed (resp. pure of weight $w$) if $O(G)$ is a sum of finite-dimensional Weil representations which are  mixed (resp. pure of weight $-w$).  Note that if $O(G)$ is pure, then it is pure of weight $0$, since the unit map $\Q_l \to O(G)$ must be Weil equivariant, so we always have a subspace of weight $0$.
\end{definition}

\begin{theorem}\label{weilred}
The natural Weil action on $\w\redpia$ is pure (of weight $0$).
\end{theorem}
\begin{proof}
Since $\w\redpia$ is reductive, its category of representations is generated under addition by the irreducible representations. Tannakian duality (\cite{tannaka}) states that $ O(\w\redpia)$ must then be dual to the pro-vector space of endomorphisms of the fibre functor from the category of representations to the category of vector spaces. Similarly, $ O(\w\redpia)\ten_{\Q_l}\bar{\Q}_l$ classifies $\bar{\Q}_l$-representations, and is dual to the fibre functor from representations over $\bar{\Q}_l$. By Schur's Lemma, scalar multiplications are the only endomorphisms of irreducible representations over $\bar{\Q}_l$. 

If we write $\End(V)$ for the space of endomorphisms of the vector space underlying $V$,
  there is then an isomorphism of $\w\redpia\by \w\redpia$-representations
$$
O(\w\redpia)\ten_{\Q_l}\bar{\Q}_l\cong \bigoplus_{V \in T} \End(V),
$$
where $T$ is the set of  all  isomorphism classes of irreducible representations $V$ of $\w\redpia$ over $\bar{\Q}_l$. By Lemma\ref{wf}, it follows that $V$ is an irreducible representation of $\pi_1(X,\bar{x})$ which is a subrepresentation of some $W(X_0,x)$-representation. This is the same as underlying a $W(X_{\bF_q^n},x)$-representation for some $n$, since $W(X_{\bF_q^n})= \pi_1(X,\bar{x})\ltimes \langle F_x^n\rangle$.

From Lafforgue's Theorem (\cite{Weil2} Conjecture 1.2.10, proved in \cite{La} Theorem VII.6 and Corollary VII.8), every  irreducible Weil representation over $\bar{\Q_l}$ is of the form
$$
V \cong P \ten \bar{\Q_l}^{(b)},
$$
for some  pure  representation $P$ of weight zero. Now,
$$
\End(V) \cong V^{\vee} \ten V \cong P^{\vee} \ten P,
$$
which is a pure $W(X_{\bF_q^n},x)$-representation of weight $0$. Therefore 
$$
\sum_{i}\End((F_x^{\sharp})^iV) =\sum_{i=0}^{n-1}\End((F_x^{\sharp})^iV) \le O(\w\redpia)\ten \bar{\Q_l})
$$
is a pure Weil subrepresentation of weight $0$. Hence $O(\w\redpia)\ten_{\Q_l}\bar{\Q}_l$ and $O(\w\redpia)$ are  also pure of weight $0$, as required.
\end{proof}

\begin{lemma}\label{weiluniv}
If $X$ is a smooth or proper variety, then $\w\algpia$ is the universal group $G$ fitting in to the diagram
$$
\algpia \to G \to \w\redpia,
$$
with $\ker(G \to \w\redpia)$ pro-unipotent.
\end{lemma}
\begin{proof}
Since $G$ and $\w\algpia$ are both quotients of $\algpia$, with $G \to \w\algpia$, it suffices to show that the composition $G \to \mathpzc{W}(X_0,x)$ is an embedding, or equivalently that the Frobenius action on $G$ is algebraic.  By Lemma \ref{algchar}, it then suffices to show that $\Hom_{\w\redpia}(V,U/[U,U])$ is finite-dimensional for all $\w\redpia$-representations $V$, where $U$ is the pro-unipotent radical of $G$. 

By studying derivations, we will see in Lemmas \ref{factorset} and \ref{radcoho}, and Proposition \ref{gpxcoho}, that   $\Hom_{\w\redpia}(U/[U,U],V)$ is just
$$
\H^1(\algpia, V)\cong\H^1(\pi_1(X,\bar{x}), V)\cong\H^1(X, \vv),
$$
  which is finite-dimensional.
\end{proof}

\begin{proposition}\label{weillevi}
The Weil action on $G$ is mixed if and only if the the induced actions on $G^{\red}$ and on the continuous dual vector space $(\Ru(G)/[\Ru(G),\Ru(G)])^{\vee}$ are mixed.
\end{proposition}
\begin{proof}
We first choose a Levi decomposition $G=G^{\red}\ltimes\Ru(G)$. The Weil action will not usually preserve this decomposition. However, for each $y \in X$, we may choose an element $u_y \in \Ru(G)$ such that $F_y':=\ad_{u_y} \circ F_y$ does preserve this Levi decomposition. The key point is that $u_y$ acts unipotently on $O(G)$.

Now, for any Weil representation $V$, the weight $a$ subrepresentation $\cW_a(V)$ of $V$ is defined as the intersection of the weight $n(y)a$ $F_y$-subrepresentations $\cW_{n(y)a}(V,F_y)$ of $V$, for all $y \in X$ and $|k(y)|=q^{n(y)}$. Since $\ad_{u_y} $ acts unipotently on $O(G)$, we deduce that  
$$
\cW_{n(y)a}(V,F_y)=\cW_{n(y)a}(V,F_y'),
$$
for all $y \in X$.

If we write $\fu$ for the (pro-nilpotent) Lie algebra of $\Ru(G)$, and let $\fu^{\vee}$ denote its continuous dual, then the isomorphism $\Ru(G)\cong \exp(\fu)$ and the Levi decomposition give us an isomorphism
$$
O(G)\cong O(G^{\red})[\fu^{\vee}]= \bigoplus_n O(G^{\red})\ten \Symm^n(\fu^{\vee}),
$$
which is $F_y'$ equivariant for all $y \in X$.

To say that a Weil representation is mixed is the same as saying that
$$
V=\bigoplus_{a\in \Z} (\bigcap_{y \in X} \cW_{n(y)a}(V,F_y)),
$$
and we have seen that for $V=O(G)$ it is equivalent to replace $F_y$ by $F_y'$. Since $O(G^{\red})$ is mixed, and this property is respected by sums and tensor operations, it suffices to show that $\fu^{\vee}$ is mixed for the $F_y'$. This is the same as being mixed for the natural action of the $F_y$ on $\fu^{\vee}$, so it suffices to show that the latter is a mixed Weil representation.   

Consider the lower central series  filtration $\Gamma_n \fu$ of $\fu$ given by 
$$
\Gamma_1 \fu:=\fu,\quad \Gamma_{n+1} \fu=[\fu,\Gamma_n \fu],
$$
so that $\fu=\Lim \fu/\Gamma_n \fu$. If $\fu_n^{\vee}:= (\fu/\Gamma_{n+1} \fu)^{\vee}$, then $\fu^{\vee}=\sum \fu_n^{\vee}$, and it only remains to show that the latter are mixed. Now there is a canonical map
$$
\fu_n^{\vee}/\fu_{n-1}^{\vee} \into \mathrm{CoLie}_n(\fu^{\vee}_1), 
$$
where $\mathrm{CoLie}_n$ is the degree $n$ homogeneous part of the free co-Lie algebra functor. Since this is a tensor operation, the right-hand side is  mixed ($\fu^{\vee}_1$ being mixed by hypothesis).

We next observe that if 
$$
0 \to V' \to V \to V''\to 0
$$
is a short exact sequence of ind-Weil representations with any two mixed, then the third is; this completes the proof.
\end{proof}

\begin{theorem}\label{mixed} If $X$ is smooth or proper, then the natural Weil action on $\w\algpia$ is mixed of non-positive weight.
\end{theorem}
\begin{proof}
By Theorem \ref{weilred} and Proposition \ref{weillevi}, it suffices to show that the Weil action on $(\Ru(\w\algpia)/[\Ru(\w\algpia),\Ru(\w\algpia)])^{\vee}$ is mixed of non-negative weight. By Lemma \ref{weiluniv} and Lemma \ref{dual}, we may alternatively describe this as 
$$
\H^1(X, \bO(\w\redpia)),
$$
where $\bO(\w\redpia)$ is the sheaf on $X$ corresponding to the vector space $O(\w\redpia)$ equipped with its left $\redpia$-action. The $\pi_1(X,\bar{x})$-action on $(\Ru(\w\algpia)/[\Ru(\w\algpia),\Ru(\w\algpia)])^{\vee}$ then comes from the right action on $O(\w\redpia)$, and  by Lemma \ref{fdual} the Frobenius action comes from the natural Frobenius action on $O(\w\redpia)$.

Now, as in Theorem \ref{weilred}, we may write
$$
O(\w\redpia)\ten_{\Q_l}\bar{\Q}_l\cong \bigoplus_{V \in T} \End(V),
$$
where $T$ is the set of  all  isomorphism classes of irreducible representations of $\w\redpia$. This is a sum of Weil representations, and  each $V$ extends to a representation of $W(X_{\bF_q^n},x)$ for some $n$, automatically compatible with the Frobenius action on $O(\w\redpia)$ (which then corresponds to the adjoint action). Since a Weil representation is pure of weight $w$ if and only if the restricted $W(X_{\bF_q^n},x)$-representation is so, it suffices to show that the $W(X_{\bF_q^n},x)$-representation
$$
\H^1(X, \vv^{\vee})\ten V
$$
is mixed for each irreducible $\pi_1(X,\bar{x})$-representation with $(F^n)^*V \cong V$. 

The group $W(X_{\bF_q^n},x)$ acts on $\H^1(X, \vv^{\vee})$ by composing the canonical map $W(X_{\bF_q^n},x) \to \Z$ with the Frobenius action arising from the Weil structure of $V$.  By Lafforgue's Theorem, we may assume that $V$ is pure of weight zero (by Schur's Lemma, note that different choices of Frobenius action on $V$ all give the same adjoint action on $\End(V)$). From Deligne's Weil II theorems (\cite{Weil2} Corollaries 3.3.4 -- 3.3.6), it then follows that $\H^1(X, \vv^{\vee})$ is mixed of non-negative weight, so $\H^1(X, \vv^{\vee})\ten V$ must also be mixed of non-negative weight, $V$ being pure of weight $0$. 
\end{proof}

\begin{corollary}
If $X$ is smooth, then the quotient map $\w\algpia \to \w\redpia$ has a unique Weil-equivariant section. 
\end{corollary}
\begin{proof}
In this case, the weights of $\H^1(X, \vv^{\vee})\ten V$ are strictly positive ($1$ or $2$), so 
$O(\w\algpia)/O(\w\redpia)$ is of strictly positive weights, giving us a decomposition
$$
O(\w\algpia)= \cW_0O(\w\algpia) \oplus \cW_{+}O(\w\algpia).
$$
Projection onto $\cW_0O(\w\algpia)=O(\w\redpia)$ yields the section.
\end{proof}

\begin{corollary}\label{pullback}
If $f:X \to Y$ is a morphism of connected varieties over $\bar{\bF}_p$, with $X$ smooth, and $\bV$ a semisimple constructible $\Q_l$-local system underlying a Weil sheaf on $Y$, then $f^{-1}\bV$ is semisimple.
\end{corollary}
\begin{proof}
If $\bV$ is of rank $n$, then it corresponds to a homomorphism $\w\varpi(Y,\bar{y})^{\red} \to \GL(n,\Q_l)$, or equivalently 
$$
O(\GL_n) \to O(\w\varpi(Y,\bar{y})^{\red})\le \cW_0O(\w\varpi(Y,\bar{y})),
$$
so $f^{-1}\bV$ must correspond to 
$$
O(\GL_n) \to \cW_0O(\w\algpia)=O(\w\redpia),
$$
as $f$ commutes with Frobenius. Therefore $f^{-1}\bV$ is semisimple.
\end{proof}

\section{Structure of the fundamental group}\label{final}

\subsection{Comparison of cohomology groups}
Fix a pro-finite group $\Gamma$, a reductive pro-algebraic group $R$ over $\Q_l$, and a Zariski-dense continuous representation $\rho\co \Gamma \to R(\Q_l)$.

We adapt the following definition from \cite{malcev} to pro-finite groups:
\begin{definition}
   Define the Malcev completion $(\Gamma)^{\rho,\mal}$   of $\Gamma$ relative to $\rho$  to be the universal diagram
$$
\Gamma \to (\Gamma,\rho)^{\mal} \xra{p} R,
$$
with $p$  a pro-unipotent extension, and the composition equal to $\rho$.
\end{definition}

\begin{remark}
Observe that if $\Gamma= \pi_1(X,\bar{x})$ and $R=\w\redpia$, with $\rho$ the canonical map, then Lemma \ref{weiluniv} shows that 
$$
(\Gamma,\rho)^{\mal}=\w\algpia.
$$
\end{remark}

\begin{lemma}\label{factorset}
For any finite-dimensional $R$-representation $V$, the canonical maps
$$
\H^i(\Gamma^{\rho,\mal},V) \to \H^i(\Gamma,V),
$$
are bijective for $i=0,1$ and injective for $i=2$. 
\end{lemma}
\begin{proof}
In both cases $\H^0(V)=V^R$ and $\H^1(V)$ is the set of continuous derivations from $\Gamma$ to $V$, which coincides with the definition of the tangent space. We now adapt the argument of \cite{haintorelli} \S 5.

Writing $G:=\Gamma^{\rho,\mal}$, we know that  $\H^2(G,V)$ is the set of isomorphism classes of extensions
$$
0 \to V \to E \to G\to 1,
$$
which pulls back to give the extension
$$
0 \to V \to E(\Q_l)\by_G\Gamma \to \Gamma\to 1
$$
of topological groups. It follows from \cite{W} 6.11.15 that $\H^2(\Gamma,V)$ classifies such extensions.

If this extension is trivial, then we have a section $\Gamma \to E(\Q_l)$, and hence a section $G \to E$, establishing injectivity. 
\end{proof}

\begin{remarks}
\begin{enumerate}
\item In the terminology of \cite{higgs}, note that the vector space-valued functors on $\Rep(R)$ given by 
$$
V \mapsto \H^i(G, V)
$$,
for a pro-unipotent extension $G \to R$, 
are the tangent space and universal obstruction space of the functor
$$
U \mapsto \Hom(G, R\ltimes U)_R/U
$$
on $\cN(R)$. 

For the latter, observe that if we have a small extension $U'\to U$ with kernel $V$, and a map $f:G \to R\ltimes \U$ over $R$, then $G \by_{ R\ltimes U}R\ltimes U'$ is an extension of $G$ by $V$, which splits if and only if  $f$ lifts to $ R\ltimes U'$.

\item 
Note that the cohomological comparison maps above can be defined in terms of derived functors, by comparing the categories of $G$-representations, $\Gamma$-representations over $\Q_l$ and $\Gamma$-representations over $\Z_l$. In particular, this means that they respect cup products.
\end{enumerate}
\end{remarks}

\begin{remark} For a pro-unipotent group $U$ equipped with an $R$-action, note that
in the terminology of \cite{higgs}, the vector space-valued functors on $\Rep(R)$ given by 
$$
V \mapsto \H^i(U, V)^R
$$
are the tangent space and universal obstruction space of $\fu$ for $i=1,2$ respectively.

For the final observation, note that if we have a small extension $\fh'\to \fh$ with kernel $V$, and a map $f:\Gamma \to R\ltimes \fh$, then $\Gamma \by_{ R\ltimes \fh}R\ltimes \fh'$ is an extension of $\Gamma$ by $V$, which splits if and only if  $f$ lifts to $\fh'$.
\end{remark}

\begin{lemma}\label{radcoho}
If $G=R\ltimes U$, for $U$ pro-unipotent, then for any $R$-representation $V$,
$$
\H^i(G,V)\cong(\H^i(U,\Q_l)\ten V)^R.
$$
\end{lemma}
\begin{proof}
The Hochschild-Serre spectral sequence gives
$$
\H^a(R, \H^b(U,V))\abuts \H^{a+b}(G,V),
$$
but $R$ is reductive, so cohomologically trivial, giving
$$
\H^i(G,V)\cong \H^i(U, V)^R\cong(\H^i(U,\Q_l)\ten V)^R,
$$
the last isomorphism following since $V$ is an $R$-representation.
\end{proof}

\begin{lemma}\label{liecoho}
If $U$ is a pro-unipotent algebraic group with associated Lie algebra $\fu$, then
$$
\H^*(U,\Q_l)\cong \H^*(\fu,\Q_l).
$$  
\end{lemma}
\begin{proof}
Observe that the categories of $U$-representations and of $\fu$-representations are equivalent.
\end{proof}

\begin{proposition}\label{gpxcoho}
Given a pointed connected algebraic variety $(Z, \bar{z})$ with \'etale fundamental group $\Gamma$, and a constructible $\Q_l$-local system $\vv$ on $Z$, the canonical maps
$$
\H^i(\Gamma, \vv_{\bar{z}}) \to \H^i(Z, \vv)
$$
are bijective for $i=0,1$, and injective for $i=2$.
\end{proposition}
\begin{proof}
It suffices to prove this for finite local systems $\ww$. If we write $\Gamma_{\alpha}$ for the finite quotients of $\Gamma$, then the canonical $\Gamma_{\alpha}$-torsors $Y_{\alpha}$ form an inverse system of varieties over $Z$, giving a Leray spectral sequence
$$
\varinjlim \H^a(\Gamma_{\alpha}, \H^b(Y_{\alpha}, \ww)) \abuts \H^{a+b}(Z, \ww).
$$
The result now follows from the observation that for $\alpha$ sufficiently large, $\H^0(Y_{\alpha}, \ww)\cong \ww_{\bar{z}}$ and $\H^1(Y_{\alpha}, \ww)=0$.
\end{proof}

\subsection{Frobenius actions}\label{frob}

We retain the conventions of \S \ref{weilact}, assuming furthermore that $X$ is either proper or smooth.

\begin{definition}
As in Theorem \ref{mixed}, $\bO(\w\redpia)$ is the sheaf of algebras on $X$ corresponding to the vector space $O(\w\redpia)$ equipped with its left $\redpia$-action. From now on, we will simply denote this sheaf by $\bO$. This is a pure Weil sheaf of weight $0$. The Frobenius actions on the  cohomology groups $\H^i(X,\bO)$ combine with the right $\redpia$-actions to make  them mixed Weil representations.
\end{definition}

\begin{theorem}\label{nfrobhull}
There is an isomorphism
$$
\Lie(\wmypia) \cong  L(\H^1(X,\bO)^{\vee})/(f(\H^2(X,\bO)^{\vee})),
$$
where $L(V)$ is the free pro-nilpotent Lie algebra on generators $V$, and  
$$
f:\H^2(X,\bO)^{\vee} \to \Gamma_2L(\H^1(X,\bO)^{\vee})
$$ 
is $R$-equivariant and preserves the (Frobenius) weight decompositions of \cite{Weil2} 3.3.7. The resulting weight decomposition on $\wmypia$ is the same as the natural Weil weight decomposition of Theorem \ref{mixed}.

Moreover, for $\ff:=L(\H^1(X,\bO)^{\vee})$,  the quotient map
$$
f:\H^2(X,\bO)^{\vee} \to \Gamma_2\ff/\Gamma_3\ff  \cong {\bigwedge}^2(\H^1(X,\bO)^{\vee})
$$
is dual to the  cup product
$$
\H^1(X,\bO) \by \H^1(X,\bO) \xra{\cup} \H^2(X,\bO).
$$

\begin{proof}
Write $G:= \w\algpia$, $R:= \w\redpia$, $U:= \wmypia$ and $\fu:= \Lie(U)$. By Lemmas \ref{radcoho}--\ref{liecoho}, Lemma \ref{dual} and Proposition \ref{gpxcoho}, we know that there is a canonical isomorphism 
$$
\H^1(\fu,\Q_l) \cong \H^1(X,\bO)
$$
of $R$-representations. Since this isomorphism is functorial, it is Frobenius-equivariant. 

We now make use of Theorem \ref{mixed}, which gives a weight decomposition on $\fu$ (which is $R$-semilinear). In fact, the theorem gives a weight decomposition on $G$, so we have an action of $\bG_m \ltimes R$ on $U$. 
To see that the Frobenius decomposition corresponds to the Weil decomposition, note that the action of $F_x \in W(X_0,x)$ determines the Weil decomposition. 

We may choose a lift of the map
$$
\fu \to \fu/[\fu,\fu]=  \H^1(\fu,\Q_l)^{\vee}\cong \H^1(X,\bO)^{\vee}
$$
as a $\bG_m \ltimes R$-representation. Writing
$\ff:=L(\H^1(X,\bO)^{\vee})$, this gives the surjection $\ff \onto \fu$.

If $J$ is the kernel of this surjection, then 
$$
J/[\ff,J] \cong \H^2(\fu,\Q_l)^{\vee},
$$
and this isomorphism is also $\bG_m \ltimes R$-linear. Once again, we may use reductivity of $\bG_m \ltimes R$ to choose a lift, giving 
$$
\H^2(\fu,\Q_l)^{\vee} \to J,
$$ 
and we define $f:\H^2(X,\bO)^{\vee} \to J$ to be the composition of this with the maps of Lemma \ref{factorset} and Proposition \ref{gpxcoho}.

Finally, the characterisation of the cup product is a standard result in Lie algebra cohomology, being dual to the map $J/[\ff,J]\to [\ff,\ff]/[\ff,[\ff,\ff]]$.
\end{proof}
\end{theorem}

\begin{corollary}\label{quadratic}
If $X$ is smooth and proper, then
$$
\Lie(\wmypia)
$$
is quadratically presented. 

In fact, there is an isomorphism of Weil representations
$$
\Lie(\wmypia) \cong  L(\H^1(X,\bO)^{\vee})/(\check{\cup}\,\H^2(X,\bO)^{\vee})),
$$
where $\check{\cup}$ is dual to the cup product.
\begin{proof}
This follows since, under these hypotheses,  \cite{Weil2} Corollaries 3.3.4--3.3.6 imply that $\H^1(X,\bO)$ is pure of weight $1$, and $\H^2(X,\bO)$ is pure of weight $2$. This makes the choices of lifts unique, and hence Frobenius-equivariant.
\end{proof}
\end{corollary}

\begin{corollary}
If $X$ is smooth and proper, there is a canonical equivalence of categories between:
\begin{enumerate}
\item the full subcategory $\C$ of the category of constructible local systems over $\Q_l$ on $X$ whose objects are subsystems of  Weil sheaves, and
\item the category of pairs $(\ww,\alpha)$, for $\ww \in \C$ semisimple and $\alpha \in \H^1(X, \End(\ww))$ with $\alpha\cup \alpha =0$.
\end{enumerate}
\end{corollary}

\begin{corollary}\label{pi1quadratic}
If $X$ is smooth and proper, then the pro-unipotent Malcev completion
$$
\pi_1(X,\bar{x})\ten\Q_l
$$
is quadratically presented. 

In fact,
$$
\cL(\pi_1(X,\bar{x}),\Q_l) \cong  L(\H^1(X,\Q_l)^{\vee})/(\check{\cup}\,\H^2(X,\Q_l)^{\vee})),
$$
where $\check{\cup}$ is dual to the cup product.
\end{corollary}
\begin{proof}
The pro-unipotent completion $\pi_1(X,\bar{x})\ten\Q_l $  is just the maximal quotient $\theta_{\sharp}\mypia$ of $\mypia$, for $\theta:\redpia \to 1$, on which $\pi_1(X,\bar{x})$ acts trivially. 
\end{proof}

\begin{example}
This implies that, for $X$ smooth and proper, the pro-$l$ quotient $\pi_1^l(X,\bar{x})$ of $\pi_1(X,\bar{x})$ cannot be the Heisenberg group 
$$
\cH_3(\Z_l)=\left\{ \begin{pmatrix} 1 & x & y\\0 & 1 & z \\ 0 & 0 & 1 \end{pmatrix} \in \GL_3(\Z_l)\right\},
$$
since this is not of quadratic presentation (in particular, this can be inferred from the non-vanishing of the Massey triple product on $\H^1(\pi_1(X,\bar{x}),\Q_l)$ ---
 see \cite{Am} Ch.3 \S 3 for criteria for a Lie algebra to be quadratically presented).
\end{example}

\begin{corollary}
If $X$ is smooth and proper, then 
$$
\cL(\pi_1(X,\bar{x}))/\Gamma_3(\cL(\pi_1(X,\bar{x}))) \ncong L(V)/\Gamma_3(L(V)),
$$ 
for any free Lie algebra $L(V)$.
\end{corollary}
\begin{proof}
As for \cite{Am} Proposition 3.25, making use of the Hard Lefschetz Theorem (\cite{Weil2} Theorem 4.1) to see that $\H^1(X,\Q_l)\by \H^1(X,\Q_l)\to \H^2(X,\Q_l)$ must be non-degenerate.
\end{proof}

\begin{corollary}\label{quartic}
If $X$ is smooth, then
$$
\Lie(\wmypia)
$$
is a quotient of the free pro-nilpotent Lie algebra $L(\H^1(X,\bO)^{\vee})$ by an ideal which is finitely generated by elements of bracket length $2,3,4$.
\begin{proof}
This follows from Theorem \ref{nfrobhull} since, under these hypotheses,  \cite{Weil2} Corollaries 3.3.4--3.3.6 imply that $\H^1(X,\bO)$ is of weights $1$ and $2$, while $\H^2(X,\bO)$ is  of weights $2, 3$ and $4$.
\end{proof}
\end{corollary}

\begin{corollary}\label{pi1quartic}
If $X$ is smooth, then
$$
\pi_1(X,\bar{x})\ten\Q_l
$$
is a quotient of the free Lie algebra $L(\H^1(X,\Q_l)^{\vee})$ by an ideal which is finitely generated by elements of bracket length $2,3,4$.
\end{corollary}

\begin{example}
Thus $\pi_1^l(X,\bar{x})$ cannot be the group
$$
\left\{ \begin{pmatrix} 1 & * & *\\0 & \ddots & * \\ 0 & 0 & 1 \end{pmatrix} \in \GL_5(\Z_l)\right\}.
$$
\end{example}

\begin{remark} If $X$ is singular and proper, weights tell us nothing about the structure of the fundamental group, since zero weights are permitted, so any equations may arise.
\end{remark} 

\subsection{Further examples}

We will now show how Theorem \ref{nfrobhull} can be used to establish stronger restrictions on the fundamental group.

\begin{corollary}\label{criterion} Let $G$ be an arbitrary  reductive $\Q_l$- algebraic group, acting on a  unipotent $\Q_l$-algebraic group $U$ defined by homogeneous equations, i.e. $\fu \cong \gr \fu$ as Lie algebras with $G$-actions. 
\begin{enumerate}
\item
If  $X$ is smooth and proper, and 
$$
\rho_2: W(X_0,x) \to (U/[U,[U,U]]) \rtimes G
$$ 
is a Zariski-dense representation, then 
$$
\rho_1: \pi_1(X,x) \to (U/[U,U]) \rtimes G
$$
 lifts to a representation 
$$
\rho: \pi_1(X,x)\to U\rtimes G.
$$

\item
If $X$ is merely smooth, and 
$$
\rho_4: W(X_0,x) \to (U/\Gamma_5 U) \rtimes G
$$ 
is a Zariski-dense representation, then 
$$
\rho_1: \pi_1(X,x) \to (U/[U,U]) \rtimes G
$$
 lifts to a representation 
$$
\rho: \pi_1(X,x)\to U\rtimes G.
$$
\end{enumerate}
\begin{proof}
As for \cite{higgs} Corollary \ref{higgs-criterion}. 
\end{proof}
\end{corollary}

\begin{remarks} 
Note that Corollaries \ref{quadratic} and \ref{quartic} imply the results  of \cite{paper1}: 

The problem considered in \cite{paper1}  is to fix a reductive representation \mbox{$\rho_0:W(X_0,x) \to G(\Q_l)$,} and consider lifts $ \rho:\pi_1(X,\bar{x}) \to G(A)$, for Artinian rings $A$. The hull of this functor is the functor
$$
A \mapsto \Hom_{\pi_1(X,\bar{x})}(\mypia, \exp(\g\ten \m_A)),
$$   
where $\g$ is the Lie algebra of $G$, regarded as the adjoint representation. It follows that this hull then has generators $\Hom_{\pi_1(X,\bar{x})}(\g, H_1)$, and relations 
$$
\Hom_{\pi_1(X,x)}(\g, H_2) \to \Symm^2 \Hom_{\pi_1(X,x)}(\g, H_1)
$$ 
given by composing the coproduct and the Lie bracket, where  
$$
H_i:=\H^i(X,\bO(\redpia))^{\vee}.
$$
\end{remarks}

\begin{definition}
A representation  $\rho:\Gamma\to G$ of a pro-finitely generated group $\Gamma$ is said to be rigid if the orbit $G(\rho) \subset \Hom(\Gamma,G)$ under the conjugation action is open in the $l$-adic topology. Observe that this is equivalent to the condition that $\H^1(\Gamma, \Lie(G))=0$, since this is the dimension of the quotient space at $[\rho]$.

A representation is properly rigid if the representation to the Zariski closure of its image is rigid.
\end{definition}

The following lemma is inspired by the observation in \cite{Simpson} that rigidity ensures that a local system on a complex projective variety is a variation of Hodge structure.
\begin{lemma}\label{rigid}
Every properly rigid representation $\rho:\pi_1(X,\bar{x})\to G$ extends to a   representation of $W(X_{k'},x)$, for some finite extension $k\subset k'$.
\begin{proof}
Replace $G$ by the Zariski closure of the image of $\rho$.
If we give the set $\N$ the multiplicative ordering, then it becomes a poset, and ${F_x^n\rho}_{n \in \N}$ is a net in $\Hom(\pi_1(X,\bar{x}),G)$. Since $F_x^n \to 1$, this net tends to $\rho$. Since $G(\rho)$ is an open neighbourhood of $\rho$, there exists an $n$ for which $F_x^n\rho \in G(\rho)$; let $F_x^n =\ad_g(\rho)$. We may now define a representation 
$$
\pi_1(X,\bar{x}) \rtimes \langle F_x^n\rangle \xra{(\rho,g)}   G,
$$
noting that the former group is $W(X_{k'},x)$, for $k\subset k'$ a degree $n$ extension.
\end{proof}
\end{lemma}

\begin{remark}
Observe that the lemma remains true under the weaker hypothesis that $\im(\H^1(\pi_1(X,\bar{x}), \Lie(\im \rho)) \to \H^1( \pi_1(X,\bar{x}), \Lie(G)))=0$.
\end{remark}

\begin{proposition}
If $X$ is smooth and proper, and $\Gamma:=\pi_1(X,\bar{x}) =\Delta \rtimes \Lambda$, let $H$ be the Zariski closure of the image of $\Lambda$ in $\Aut(\Delta\ten \Q_l)$. If $H$ is reductive, $\H^1(\L, \Lie(H))=0$, and $\Hom_{\L}(\Delta/[\Delta,\Delta], \Lie(H))=0$,    then $\Delta\ten \Q_l$ is quadratically presented.
\begin{proof} 
 First observe that the representation $\rho:\Gamma \to \Lambda \to H$ is rigid. This follows because the condition  $\H^1(\L, \Lie(H))=0$ ensures that $\L \to H$ is rigid, so any representation $\Gamma \to H \ltimes \eps\Lie(H)$ (for $\eps^2=0$) must be conjugate to one which restricts to $\rho$ on $\L$. The image of $\Delta$ must also lie in $\Lie(H)$, so the representation is determined  by an element of $\Hom_{\L}(\Delta/[\Delta,\Delta], \Lie(H))=0$, so it must be $\rho$. Therefore, by Lemma \ref{rigid}, $\rho$ extends to a Weil representation (possibly after changing the base field). 

Hence $\rho$ factors as $\pi_1(X,\bar{x}) \to\w\redpia \xra{\theta} H$.
 Now, by  Corollary \ref{quadratic}, we know that
$$
\theta_{\sharp}Lie(\wmypia)
$$
must be quadratically presented. The proof now proceeds as in \cite{higgs} Proposition \ref{higgs-crit}.
\end{proof}
\end{proposition}

\begin{example} Let $\mathfrak{d}$ be the free $\hat{\Z}$-module
$$
\mathfrak{d}:=\hat{\Z}x \oplus  \hat{\Z}y \oplus \half \hat{\Z}[x,y],
$$
which has the structure of a Lie algebra, with $[x,y]$ in the centre.
 The Campbell-Baker-Hausdorff formula enables us to regard $\Delta:=\exp(\mathfrak{d})$ as the profinite group with underlying set $\mathfrak{d}$ and product
$$
a\cdot b=a + b +\half [a,b],
$$
since all higher brackets vanish.

Let $\exp(\fh):=\Delta\ten \Q_l$; this is isomorphic to the  three-dimensional $l$-adic Heisenberg group.  

Observe that $\Aut(\Delta \ten \Q_l) \cong \GL_2(\Q_l)$, and that $\SL_2(\hat{\Z})$ acts on $\Delta$ by the formula:
$$
A(v,w) := (Av, (\det A) w)=(Av,w),
$$
for $v \in \hat{\Z}x \oplus \hat{\Z}y$ and $w \in \half \hat{\Z}[x,y]$. 

The group
 $$
\Gamma:=\Delta \rtimes \SL_2(\hat{\Z});
$$
cannot be the geometric fundamental group of any smooth proper variety defined over the algebraic closure of a finite field.
\begin{proof}
We wish to show that $\Gamma \to \Aut(\Delta \ten \Q_l)$ is properly rigid. For this, it will suffice to show that $\L \to  \Aut(\Delta \ten \Q_l)$ is properly rigid, and that $\Hom_{\L}(\Delta/[\Delta,\Delta], \mathfrak{sl}_2(\Q_l))=0$. 

To prove the first, observe that 
$$
\SL_2(\hat{\Z})=\prod_{\nu \text{ prime}}\SL_2(\Z_{\nu}),
$$
and that only pro-$l$ groups contribute to cohomology. We need to show that the only derivations $\SL_2(\Z_l) \to \mathfrak{sl}_2(\Q_l)$ are inner derivations.  Now, for $N$ sufficiently large, $\exp:l^N\mathfrak{sl}_2(\Z_l) \to  \SL_2(\Z_l)$ converges, and it follows from the simplicity of  $\mathfrak{sl}_2(\Z_l)$ that any derivation must agree with an inner derivation when restricted to $\exp(l^N\mathfrak{sl}_2(\Z_l))$. Since this is a subgroup of finite index, and $\mathfrak{sl}_2(\Q_l)$ is torsion-free, the  derivation and inner derivation must agree on the whole of  $\SL_2(\Z_l)$, as required. 

To prove the second, observe that $\Q_l^2$ and $\mathfrak{sl}_2(\Q_l)$ are distinct irreducible $\SL_2(\Z_l)$-representations.

We therefore conclude from the previous proposition that $\Gamma$
cannot be the fundamental group of any smooth proper variety defined over the algebraic closure of a finite field, since  the action of $\SL_2(\hat{\Z})$ on $\fh$ is semisimple, hence reductive, and $\fh$ is not quadratically presented. 

Alternatively, we could use Corollary \ref{criterion} to prove that $\Gamma$ is not such a group. Let $G=\SL_2(\Q_l)$, $\fu=L(\Q_l^2)$ and $U=\exp(\fu)$. Observe that $\fh \cong \fu/[\fu,[\fu,\fu]]$, and let $\rho_2$ be the standard embedding
$$
\rho_2: \Delta \rtimes \SL_2(\hat{\Z}) \to \exp(\fh) \rtimes \SL_2(\Q_l),
$$
which extends to a Weil representation by Lemma \ref{rigid} and the above calculation. 

Since all triple commutators vanish in $H$, this does not lift to a representation
$$
\rho: \Delta \rtimes \SL_2(\hat{\Z}) \to U \rtimes \SL_2(\Q_l).
$$
\end{proof}

Note that Corollary \ref{pi1quartic} cannot be used to exclude this group ---  the abelianisation of $\Gamma$ is a torsion group, as $\SL_2$ acts irreducibly on the abelianisation of $\fh$, so $\Gamma\ten \Q_l =1$, which is quadratically presented.
\end{example}

\begin{example}
Let $\mathfrak{d}$ be  the free $\hat{\Z}$-module on generators 
$$
x,y, \frac{1}{2}[xy], \frac{1}{12}[x[xy]], \frac{1}{12}[y[xy]], \frac{1}{24}[x[x[xy]]],\frac{1}{24}[x[y[xy]]],\frac{1}{24}[y[y[xy]]];
$$
this has a Lie algebra structure, with all quintuple commutators vanishing. We define $\Delta:=\exp(\mathfrak{d})$, the group whose underlying set is $\mathfrak{d}$, given a group structure via the truncated Campbell-Baker-Hausdorff formula:
$$
a\cdot b=a+b+\frac{1}{2}[a,b]+\frac{1}{12}([a,[a,b]]-[b,[a,b]]) -\frac{1}{24}[a,[b,[a,b]]].
$$
Again, let $\L:=\SL_2(\hat{\Z})$, acting in the natural way on $\hat{\Z}x \oplus \hat{\Z}y$, with the action extending to the whole of $\Delta$ via the laws of Lie algebras. Then 
$$
\Gamma:=\Delta \rtimes \SL_2(\hat{\Z})
$$
cannot be the geometric fundamental group of any smooth variety defined over the algebraic closure of a finite field.
\begin{proof}
Let $G=\SL_2(\Q_l)$, and let  $\fu=L(\Q_l^2)$ and $U=\exp(\fu)$. Observe that 
$$
\fh:=\cL(\Delta, \Q_l) \cong \fu/\Gamma_5\fu,
$$ 
and let $\rho_4$ be the standard embedding
$$
\rho_4: \Delta \rtimes \SL_2(\hat{\Z}) \to \exp(\fh) \rtimes \SL_2(\Q_l),
$$
which extends to a Weil representation by Lemma \ref{rigid} and  calculation in the previous example. 

Since all quintuple commutators vanish in $H$, this does not lift to a representation
$$
\rho: \Delta \rtimes \SL_2(\hat{\Z}) \to U \rtimes \SL_2(\Q_l),
$$
which gives a contradiction, by Corollary \ref{criterion}.
\end{proof}
\end{example}

\bibliographystyle{alphanum}
\addcontentsline{toc}{section}{Bibliography}
\bibliography{references.bib}
\end{document}

%% file: newcommands.tex
\newcommand\ten{\otimes}

\newcommand\eps{\epsilon}

\newcommand\Ru{\mathrm{R_u}}

\renewcommand\H{\mathrm{H}}

\newcommand\N{\mathbb{N}}
\newcommand\Z{\mathbb{Z}}
\newcommand\Q{\mathbb{Q}}

\newcommand\vv{\mathbb{V}}
\newcommand\ww{\mathbb{W}}

\newcommand\bF{\mathbb{F}}
\newcommand\bG{\mathbb{G}}

\newcommand\bO{\mathbb{O}}

\newcommand\bV{\mathbb{V}}

\newcommand\C{\mathcal{C}}
\newcommand\U{\mathcal{U}}

\newcommand\cH{\mathcal{H}}

\newcommand\cL{\mathcal{L}}

\newcommand\cN{\mathcal{N}}

\newcommand\cW{\mathcal{W}}

\renewcommand\L{\Lambda}

\newcommand\m{\mathfrak{m}}

\newcommand\g{\mathfrak{g}}

\newcommand\ff{\mathfrak{f}}

\newcommand\fh{\mathfrak{h}}

\newcommand\fu{\mathfrak{u}}

\newcommand\Hom{\mathrm{Hom}}

\newcommand\End{\mathrm{End}}

\newcommand\Aut{\mathrm{Aut}}

\newcommand\im{\mathrm{Im\,}}

\newcommand\Ind{\mathrm{Ind}}

\newcommand\mal{\mathrm{Mal}}

\newcommand\Spec{\mathrm{Spec}\,}

\newcommand\ad{\mathrm{ad}}

\newcommand\Lim{\varprojlim}

\newcommand\into{\hookrightarrow}
\newcommand\onto{\twoheadrightarrow}
\newcommand\abuts{\implies}
\newcommand\xra{\xrightarrow}

\newcommand\alg{\mathrm{alg}}

\newcommand\by{\times}

\newcommand\Rep{\mathrm{Rep}}

\newcommand\Symm{\mathrm{Symm}}
\newcommand\SL{\mathrm{SL}}
\newcommand\GL{\mathrm{GL}}

\newcommand\half{\frac{1}{2}}

\newcommand\gr{\mathrm{gr}}

\newcommand\redpia{\varpi_1^{\mathrm{red}}(X,\bar{x})}
\newcommand\algpia{\varpi_1(X,\bar{x})}
\newcommand\mypia{R_u(\varpi_1(X,\bar{x}))}
\newcommand\w{{}^W\!}
\newcommand\wmypia{R_u(\w\varpi_1(X,\bar{x}))}
\renewcommand\alg{\mathrm{alg}}
\newcommand\red{\mathrm{red}}

\newcommand\Lie{\mathrm{Lie}}

\newcommand\co{\colon\thinspace}

%% file: newthms.tex
\newtheorem{theorem}{Theorem}[section]
\newtheorem{proposition}[theorem]{Proposition}
\newtheorem{corollary}[theorem]{Corollary}
\newtheorem{lemma}[theorem]{Lemma}
\newtheorem{theorem*}{Theorem}
\newtheorem{proposition*}[theorem*]{Proposition}
\newtheorem{corollary*}[theorem*]{Corollary}
\newtheorem{lemma*}[theorem*]{Lemma}

\theoremstyle{definition}
\newtheorem{definition}[theorem]{Definition}

\newtheorem{definition*}[theorem*]{Definition}

\theoremstyle{remark}
\newtheorem{example}[theorem]{Example}

\newtheorem{remark}[theorem]{Remark}
\newtheorem{remarks}[theorem]{Remarks}

\newtheorem{example*}[theorem*]{Example}
\newtheorem{examples*}[theorem*]{Examples}
\newtheorem{remark*}[theorem*]{Remark}
\newtheorem{remarks*}[theorem*]{Remarks}
\newtheorem{exercise*}[theorem*]{Exercise}